\documentclass[11pt]{amsart}
\usepackage{amssymb}
\usepackage[all]{xy}

\newcommand{\myendbibitem}{\relax}

\usepackage[linktocpage=true,hypertex]{hyperref}
\swapnumbers
\newtheorem{thm}[equation]{Theorem}
\newtheorem{prop}[equation]{Proposition}
\newtheorem{lem}[equation]{Lemma}
\newtheorem{cor}[equation]{Corollary}

\theoremstyle{definition}
\newtheorem{defn}[equation]{Definition}
\newtheorem{remark}[equation]{Remark}
\newtheorem{example}[equation]{Example}

\numberwithin{equation}{section}

\newcommand{\bbA}{{\mathbb A}}
\newcommand{\bbP}{{\mathbb P}}
\newcommand{\calI}{{\mathcal I}}
\newcommand{\calZ}{{\mathcal Z}}
\newcommand{\id}{\operatorname{id}}
\newcommand{\inv}{^{-1}}
\newcommand{\Ker}{\operatorname{Ker}}
\newcommand{\lra}{\longrightarrow}
\newcommand{\PIdeg}{\operatorname{PIdeg}}
\newcommand{\BS}{\operatorname{BS}}
\newcommand{\RMaps}{\operatorname{RMaps}}

\newcommand{\Spec}{\operatorname{Spec}}
\newcommand{\tr}{\operatorname{tr}}

\newcommand{\Cmn}{C_{m,n}}
\newcommand{\Gmn}{G_{m,n}}
\newcommand{\Qmn}{Q_{m,n}}
\newcommand{\Tmn}{T_{m,n}}
\newcommand{\Umn}{U_{m,n}}

\newcommand{\M}{\operatorname{M}}
\newcommand{\Mat}{{\operatorname{M}}}
\newcommand{\Mn}{\Mat_n}
\newcommand{\Mnl}{(\Mn)^l}
\newcommand{\Mnm}{(\Mn)^m}
\newcommand{\GL}{{\operatorname{GL}}}
\newcommand{\PGL}{{\operatorname{PGL}}}
\newcommand{\PGLn}{{\operatorname{PGL}_n}}

\newcommand{\Var}{\textit{Var}}
\newcommand{\PI}{\textit{PI}}
\newcommand{\Bir}{\textit{Bir}}
\newcommand{\CS}{\textit{CS}}

\newcommand{\numberedpar}{%%%
        \par\refstepcounter{equation}\medskip\noindent{\bf\theequation. }}

\begin{document}

\title[Polynomial identity rings]%
{Polynomial identity rings as rings of functions}

\date{August 5, 2005}

\author{Z. Reichstein}
\address{Department of Mathematics, University of British Columbia,
       Vancouver, BC V6T 1Z2, Canada}
\email{reichst@math.ubc.ca}
\urladdr{www.math.ubc.ca/$\stackrel{\sim}{\phantom{.}}$reichst}
\thanks{Z. Reichstein was supported in part by an NSERC research
  grant}

\author{N. Vonessen} \address{Department of Mathematical Sciences,
  University of Montana, Missoula, MT 59812-0864, USA}
\email{Nikolaus.Vonessen@umontana.edu}
\urladdr{www.math.umt.edu/vonessen}
\thanks{N.\ Vonessen gratefully acknowledges the support of the
  University of Montana and the hospitality of the University of
  British Columbia during his sabbatical in 2002/2003, when part of
  this research was done.}

\begin{abstract}
  We generalize the usual relationship between irreducible Zariski
  closed subsets of the affine space, their defining ideals,
  coordinate rings, and function fields, to a non-commutative setting,
  where ``varieties" carry a $\PGLn$-action, regular and rational
  ``functions" on them are matrix-valued, ``coordinate rings" are
  prime polynomial identity algebras, and ``function fields" are
  central simple algebras of degree~$n$.  In particular, a prime
  polynomial identity algebra of degree $n$ is finitely generated if
  and only if it arises as the ``coordinate ring" of a ``variety'' in this
  setting.  For $n = 1$ our definitions and results reduce to those of
  classical affine algebraic geometry.
\end{abstract}


%%%%%%%%%%%%%%%%%%%%%%%%%%%%%%%%%%%%%%%%%%%%%%%%%%%%%%%%%%%%%%%%%%%%%%
% 14A10 (1973-now) Varieties and morphisms
% 14L30 Group actions on varieties and schemes
% 16R20 (1991-now) Semiprime p.i. rings, rings embeddable
% in matrices over commutative rings
% 16R30 (1991-now) Trace rings and invariant theory
%%%%%%%%%%%%%%%%%%%%%%%%%%%%%%%%%%%%%%%%%%%%%%%%%%%%%%%%%%%%%%%%%%%%%%


\maketitle

\tableofcontents

\section{Introduction}

Polynomial identity rings (or PI-rings, for short)
are often viewed as being ``close to commutative"; 
they have large centers, and their structure (and 
in particular, their maximal spectra) have been successfully
studied by geometric means (see the references at the end 
of the section).  In this paper we revisit this 
subject from the point of view of classical affine 
algebraic geometry.  We will show that the usual 
relationship between irreducible Zariski closed subsets
of the affine space, their defining ideals, coordinate
rings, and function fields, can be extended to the setting
of PI-rings.

Before proceeding with the statements of our main results,
we will briefly introduce the objects that will play the roles
of varieties, defining ideals, coordinate rings, etc.
Throughout this paper we will work over
an algebraically closed base field $k$ of
characteristic zero. We also fix an integer $n \ge 1$, which
will be the PI-degree of most of the rings we will consider.
We will write $\Mn$ for the matrix algebra $\Mn(k)$.
The vector space of $m$-tuples of $n \times n$-matrices will
be denoted by $\Mnm$; we will always assume that $m \ge 2$.
The group $\PGLn$ acts on $\Mnm$ by simultaneous conjugation.
The $\PGLn$-invariant dense open subset
\[ \text{$U_{m, n} = \{ (a_1, \dots, a_m) \in \Mnm \, |
\, a_1, \dots, a_m$ generate $\Mn$ as $k$-algebra$\}$ }  \]
of $\Mnm$ will play the role of the affine space $\bbA^m$ in
the sequel. (Note that $U_{m, 1} = \bbA^m$.)
The role of affine algebraic varieties will be played by
$\PGLn$-invariant closed subvarieties of $U_{m, n}$;
for lack of a better term, we shall call such objects
$n$-{\em varieties}; see Section~\ref{sect.defn}. (Note that,
in general, $n$-varieties are not affine in the usual sense.)
The role of the polynomial ring $k[x_1, \dots, x_m]$ will be played
by the algebra $\Gmn = k \{ X_1, \dots, X_m \}$ of $m$ generic
$n\times n$ matrices, see~\ref{prel.Gmn}.
Elements of $\Gmn$ may be thought of as $\PGLn$-equivariant maps
$\Mnm \lra \Mn$; if $n = 1$ these are simply the polynomial maps
$k^m \lra k$.  Using these maps, we define, in a manner analogous
to the commutative case, the associated ideal $\calI(X) \subset \Gmn$,
the PI-coordinate ring $k_n[X] = G_{m, n}/\calI(X)$, and
the central simple algebra $k_n(X)$ of rational functions
on an irreducible $n$-variety $X \subset U_{m, n}$; see
Definitions~\ref{def.n-var1} and~\ref{def.ratl.1}.
We show that Hilbert's Nullstellensatz continues to
hold in this context; see Section~\ref{sect.nullstellensatz}.
We also define the notions of a regular map (and,
in particular, an isomorphism) $X \lra Y$ and a rational
map (and, in particular, a birational isomorphism)
$X \dasharrow Y$ between $n$-varieties
$X \subset U_{m, n}$ and $Y \subset U_{l, n}$; see
Definitions~\ref{def.n-var2} and~\ref{def.ratl.2}.

In categorical language our main results can be
summarized as follows. Let

\smallskip $\Var_n$ be the category of irreducible $n$-varieties, with
regular maps of $n$-varieties as morphisms (see
Definition~\ref{def.n-var2}), and

\smallskip $\PI_n$ be the category of finitely generated prime
$k$-algebras of PI-degree $n$ (here the morphisms are the usual
$k$-algebra homomorphisms).

\begin{thm} \label{thm1}
The functor defined by
\[ \begin{array}{rcl}  X &\mapsto &k_n[X] \\
     (f \colon X \lra Y) &\mapsto &(f^* \colon k_n[Y] \lra k_n[X])
\end{array} \]
is a contravariant equivalence of categories
between $\Var_n$ and $\PI_n$.
\end{thm}

\smallskip
In particular, every finitely generated prime PI-algebra is
the coordinate ring of a uniquely determined $n$-variety; see
Theorem~\ref{thm.n-Var}. For a proof of Theorem~\ref{thm1},
see Section~\ref{sect.reg}.

\smallskip
Every $n$-variety is, by definition, an algebraic variety with
a generically free $\PGLn$-action. It turns out that, up to
birational isomorphism, the converse holds as well;
see Lemma~\ref{lem.rat3}. To summarize our results
in the birational context, let

\smallskip
$\Bir_n$ be the category of irreducible generically free
$\PGLn$-varieties, with dominant rational $\PGLn$-equivariant maps as
morphisms, and

\smallskip
$\CS_n$ be the category of central simple
algebras $A$ of degree $n$, such that the center of $A$
is a finitely generated field extension of $k$.
Morphisms in $\CS_n$ are $k$-algebra homomorphisms
(these are necessarily injective).

\begin{thm} \label{thm2}
The functor defined by
\[ \begin{array}{rcl}      X &\mapsto &k_n(X) \\
   (g \colon X \dasharrow Y) &\mapsto &(g^* \colon k_n(Y)
         \hookrightarrow k_n(X)) \end{array} \]
is a contravariant equivalence of categories between $\Bir_n$ and $\CS_n$.
\end{thm}

Here for any $\PGLn$-variety $X$, $k_n(X)$ denotes the $k$-algebra of
$\PGLn$-equivariant rational maps $X \dasharrow \Mn$ (with addition
and multiplication induced from $\Mn$), see
Definition~\ref{def.2.kn(X)}.  If $X$ is an irreducible $n$-variety then
$k_n(X)$ is the total ring of fractions of
$k_n[X]$, as in the commutative case; see Definition~\ref{def.ratl.1} and
Proposition~\ref{prop.rat2}.  Note also that $g^*(\alpha)$ stands for
$\alpha \circ g$ (again, as in the commutative case).
For a proof of Theorem~\ref{thm2}, see Section~\ref{sect.gener}.

Note that for $n = 1$, Theorems~\ref{thm1} and~\ref{thm2} are classical
results of affine algebraic geometry; cf.,
e.g.,~\cite[Corollary 3.8 and Theorem 4.4]{hartshorne}.

It is well known that
central simple algebras $A/K$ of degree $n$ are in a natural bijection
with $n-1$-dimensional Brauer-Severi varieties over $K$ and
(if $K/k$ is a finitely generated field extension) 
with generically free 
$\PGLn$-varieties $X/k$ such that $k(X)^{\PGLn} \simeq K$.
Indeed, all three are parametrized by the Galois cohomology
set $H^1(K, \PGLn)$; for details, see Section~\ref{sect.bs}.
Theorem~\ref{thm2} may thus be viewed as a way of explicitly
identifying generically free $\PGLn$-varieties $X$ with central simple
algebras $A$, 
without going through $H^1(K, \PGLn)$. The fact that
the map $X \mapsto k_n(X)$ is bijective was proved 
in \cite[Proposition 8.6 and Lemma 9.1]{reichstein4};
here we give a more conceptual proof and show that
this map is, in fact, a contravariant functor.
In Section~\ref{sect.bs} we show how to construct the
Brauer-Severi variety of $A$ directly from $X$.

Many of the main themes of this paper
(such as the use of $\PGLn$-actions, generic matrices, 
trace rings, and affine geometry in the study of polynomial 
identity algebras)
were first systematically explored in the pioneering 
work of Amitsur, Artin and 
Procesi~\cite{amitsur:procesi, artin1, procesi2, procesi3, procesi1}
in the 1960s and 70s.  Our approach here was influenced by these, 
as well as other papers in this area, such 
as~\cite{amitsur:small, artin2, as1, razmyslov}.  
In particular, Proposition~\ref{prop.strong-null} is similar 
in spirit to Amitsur's Nullstellensatz~\cite{amitsur}
(cf.\ Remark~\ref{rem.Amitsur}), and
Theorem~\ref{thm1} to Procesi's functorial description
of algebras satisfying the $n$th Cayley-Hamilton 
identity~\cite[Theorem 2.6]{procesi4}. (For a more geometric
statement of Procesi's theorem, along the lines of
our Theorem~\ref{thm1}, see~\cite[Theorem 2.3]{lebruyn}.)
We thank L. Small and the referee for bringing some 
of these connections to our attention.

\section{Preliminaries}
\addtocontents{toc}{\ }
\addtocontents{toc}{\hbox{}\ \ \quad\ref{prel.mi} Matrix Invariants
   \hfill\pageref{prel.mi}\\}
\addtocontents{toc}{\hbox{}\ \indent\ \ \quad\ref{prel.Gmn} The ring
  of generic matrices and its trace ring\hfill\pageref{prel.Gmn}\\}
\addtocontents{toc}{\hbox{}\ \indent\ \ \quad\ref{prel.cp} Central
  polynomials\hfill\pageref{prel.cp}}

In this section we review some known results about
matrix invariants and related PI-theory.

\numberedpar\label{prel.mi}{\bf Matrix invariants.}
Consider the diagonal action of $\PGLn$ on the space
$\Mnm$ of $m$-tuples of $n \times n$-matrices. We
shall denote the ring of invariants for this action by
$C_{m, n} = k[\Mnm]^{\PGLn}$, and the affine variety
$\Spec(C_{m, n})$ by $Q_{m, n}$.
It is known that $C_{m, n}$ is generated as a $k$-algebra, by elements
of the form $(A_1, \dots, A_m) \mapsto \tr(M)$, where $M$ is
a monomial in $A_1, \dots, A_m$ (see~\cite{procesi1}); however, we shall
not need this fact in the sequel.  The inclusion $C_{m, n} \hookrightarrow
k[\Mnm]$ of $k$-algebras induces the categorical quotient map
\begin{equation} \label{e.catq}
\pi \colon \Mnm \lra Q_{m, n} \, .
\end{equation}
We shall need the following facts about this map in the sequel.
Recall from the introduction that we always assume $m \ge 2$,
the base field $k$ is algebraically closed and of characteristic zero,
and
\[ \text{$U_{m, n} = \{ (a_1, \dots, a_m) \in \Mnm \, |
\, a_1, \dots, a_m$ generate $\Mn$ as $k$-algebra $\}$.}  \]

\begin{prop} \label{prel.matr-inv}
\begin{itemize}
\item[]\hspace{-3em}{\upshape(a)} If $x \in U_{m, n}$ then $\pi^{-1}(\pi(x))$ is the
  $\PGLn$-orbit of $x$.
\item[\upshape(b)]$\PGLn$-orbits in $\Umn$ are closed in $\Mnm$.
\item[\upshape(c)]$\pi$ maps closed $\PGLn$-invariant sets in $\Mnm$
  to closed sets in $Q_{m, n}$.
\item[\upshape(d)]$\pi(U_{m,n})$ is Zariski open in $\Qmn$.
\item[\upshape(e)]If $Y$ is a closed irreducible subvariety of $Q_{m,
    n}$ then $\pi^{-1}(Y) \cap U_{m, n}$ is irreducible in $\Mnm$.
\end{itemize}
\end{prop}

\begin{proof} (a) is proved in~\cite[(12.6)]{artin1}.

(b) is an immediate consequence of (a).

(c) is a special case of~\cite[Corollary to Theorem 4.6]{pv}.

(d) It is easy to see that $U_{m,n}$ is Zariski open in $\Mnm$.
Let $U_{m, n}^c$ be its complement in $\Mnm$. By (c),
$\pi(U_{n, m}^c)$ is closed in $Q_{m, n}$ and by (a),
$\pi(U_{m, n}) = Q_{m,n} \setminus \pi(U_{m, n}^c)$.

(e) Let $V_1, \dots, V_r$ be the irreducible components of $\pi^{-1}(Y)$
in $\Mnm$.  Since $\PGLn$ is connected, each $V_i$ is
$\PGLn$-invariant.  By part (c), $\pi(V_1)$, \dots, $\pi(V_r)$ are closed
subvarieties of $Q_{m, n}$ covering $Y$.  Since $Y$ is irreducible,
we may assume, after possibly renumbering
$V_1, \dots, V_r$, that $Y = \pi(V_1)$. It suffices to show that
\[ \pi^{-1}(Y) \cap U_{m, n} = V_1 \cap U_{m, n} \, ; \]
since $V_1 \cap U_{m, n}$ is irreducible,
this will complete the proof of part~(e).
As $V_1$ is an irreducible component of $\pi^{-1}(Y)$, we clearly have
\[ V_1 \cap U_{m, n}  \subset \pi^{-1}(Y) \cap U_{m, n} \, . \]
To prove the opposite inclusion,
let $y \in \pi^{-1}(Y) \cap U_{m, n}$. We want to show $y \in V_1$.
Since $\pi(y) \in Y = \pi(V_1)$, there is a point $v \in V_1$ such that
$\pi(y) = \pi(v)$. That is, $v$ lies in $\pi^{-1}(\pi(y))$, which,
by part (a), is the $\PGLn$-orbit of $y$. In other words,
$y = g \cdot v$ for some $g \in \PGLn$. Since $V_1$ is $\PGLn$-invariant,
this shows that $y \in V_1$, as claimed.
\end{proof}

\numberedpar\label{prel.Gmn}{\bf The ring of generic matrices and its
  trace ring.}
Consider $m$ generic matrices
\[ X_1 = (x_{ij}^{(1)})_{i, j = 1, \dots, n}\,, \dots,\,
X_m = (x_{ij}^{(m)})_{i, j = 1, \dots, n} \, , \]
where $x_{ij}^{(h)}$ are $m n^2$ independent variables over
the base field $k$. The $k$-subalgebra generated by $X_1, \dots, X_m$
inside $\Mn(k[x_{ij}^{(h)}])$ is called {\em the algebra of $m$
generic $n \times n$-matrices} and is denoted by $\Gmn$. If
the values of $n$ and $m$ are clear from the context, we will simply
refer to $\Gmn$ as the algebra of generic matrices.

The trace ring of $\Gmn$ is denoted by $T_{m, n}$; it is
the $k$-algebra generated, inside $\Mn(k[x_{ij}^{(h)}])$
by elements of $\Gmn$ and their traces.
Elements of $\Mn(k[x_{ij}^{(h)}])$ can be naturally viewed as
regular (i.e., polynomial) maps $\Mnm \lra \Mn$. (Note that
$k[x_{ij}^{(h)}]$ in the coordinate ring of $\Mnm$.)
Here $\PGLn$ acts on both $\Mnm$ and $\Mn$ by simultaneous conjugation;
Procesi~\cite[Section 1.2]{procesi1} noticed that
$T_{m, n}$ consists precisely of those maps
$\Mnm \lra \Mn$ that are equivariant with respect to this action.
(In particular, the $i$-th generic matrix $X_i$ is
the projection to the $i$-th component.) In this way
the invariant ring $C_{m, n} = k[\Mnm]^{\PGLn}$ which we considered
in~\ref{prel.mi}, is naturally identified with
the center of $T_{m, n}$ via $f \mapsto f I_{n \times n}$.

We now recall the following definitions.

\begin{defn} \label{prel.def1}
(a) A prime PI-ring is said to have PI-degree $n$ if its total
ring of fractions is a central simple algebra of degree~$n$.

\smallskip
(b) Given a ring $R$, $\Spec_n(R)$ is defined as the set of prime
ideals $J$ of $R$ such that $R/J$ has PI-degree $n$; cf.  e.g.,
\cite[p.~58]{procesi2} or \cite[p.~75]{rowen:PI}.
\end{defn}

The following lemma shows that $\Spec_n(\Gmn)$ and $\Spec_n(\Tmn)$
are closely related.

\begin{lem}\label{lem:TraceRingofCoordRing:new}
  The assignment $J \mapsto J \cap G_{m, n}$ defines a bijective
  correspondence between $\Spec_n(T_{m, n})$ and $\Spec_n(G_{m, n})$.
  In addition, for any prime ideal $J\in\Spec_n(\Tmn)$, we have the
  following: 
\begin{itemize}
\item[\upshape(a)]The natural projection $\phi \colon T_{m, n} \lra
  T_{m, n}/J$ is trace-preserving, and $\Tmn/J$ is the trace ring of
  $\Gmn/(J \cap G_{m, n})$.
\item[\upshape(b)]$\tr(p) \in J$ for every $p \in J$.
\end{itemize}
\end{lem}

\begin{proof}
  The first assertion and part (a) are special cases of results proved
  in \cite[\S2]{amitsur:small}. Part~(b) follows from~(a), since for
  any $p\in J$, $\phi(\tr(p)) = \tr(\phi(p)) = \tr(0) = 0$.  In other
  words, $\tr(p) \in \Ker(\phi) = J$, as claimed.
\end{proof}

\numberedpar\label{prel.cp}{\bf Central polynomials.}  
We need to construct central polynomials with certain non-vanishing
properties.  We begin by recalling two well-known facts from the
theory of rings satisfying polynomial identities.

\begin{prop}\label{prel.cp.prop}
  {\upshape(a)} Let $k\{x_1, \ldots, x_m \}$ be the free associative
  algebra.  Consider the natural homomorphism $k\{x_1, \dots, x_m \}
  \lra \Gmn$, taking $x_i$ to the $i$-th generic matrix $X_i$. The
  kernel of this homomorphism is precisely the ideal of polynomial
  identities of $n \times n$-matrices in~$m$ variables.
  
  \smallskip {\upshape(b)} Since $k$ is an infinite field, all prime
  $k$-algebras of the same PI-degree satisfy the same polynomial
  identities {\upshape(}with coefficients in $k${\upshape)}.
\end{prop}

\begin{proof}
See~\cite[pp.~20-21]{procesi2} or \cite[p.~16]{rowen:PI}
for a proof of part (a) and \cite[pp.~106-107]{rowen:ringII}
for a proof of part (b). 
\end{proof}

For the convenience of the reader and lack of a suitable reference, we
include the following definition.

\begin{defn}
  An ($m$-variable) {\em central polynomial for $n\times n$ matrices}
  is an element $p=p(x_1, \dots, x_m) \in k\{x_1, \ldots, x_m \}$
  satisfying one of the following equivalent conditions:
  \begin{itemize}
  \item[(a)]$p$ is a polynomial identity of $\M_{n-1}$, and the
    evaluations of $p$ in $\Mn$ are central (i.e., scalar matrices)
    but not identically zero.
  \item[(b)]$p$ is a polynomial identity for all prime $k$-algebras of
    PI-degree~$n-1$, and the evaluations of $p$ in every prime
    $k$-algebra of PI-degree~$n$ are central but not identically zero.
  \item[(c)]$p$ is a polynomial identity for all prime $k$-algebras of
    PI-degree~$n-1$, and the canonical image of $p$ in $\Gmn$ is a
    nonzero central element.
  \item[(d)]The constant coefficient of $p$ is zero, and the canonical
    image of $p$ in $\Gmn$ is a nonzero central element.
  \end{itemize}
\end{defn}

That the evaluations of $p$ in an algebra $A$ are central is
equivalent to saying that $x_{m+1}p-px_{m+1}$ is a polynomial identity
for $A$, where $x_{m+1}$ is another free variable.  Thus the
equivalence of (a)---(c) easily follows from
Proposition~\ref{prel.cp.prop}.  The equivalence of (c) and (d)
follows from \cite[p.~172]{procesi2}.  The existence of central
polynomials for $n \times n$-matrices was established independently by
Formanek and Razmyslov; see~\cite{formanek:cbms}.  Because of
Proposition~\ref{prel.cp.prop}(a), one can think of $m$-variable
central polynomials of $n\times n$ matrices as nonzero central
elements of $\Gmn$ (with zero constant coefficient). 

The following lemma, establishing the existence of central polynomials
with certain non-vanishing properties, will be repeatedly used in the
sequel.

\begin{lem} \label{prel3.1}
  Let $A_1, \dots, A_r \in U_{m, n}$. Then there exists a central
  polynomial $s=s(X_1, \dots, X_m)\in\Gmn$ for $n \times
  n$-matrices such that $s(A_i) \ne 0$ for $i = 1, \dots, r$.
  In other words, each $s(A_i)$ is a non-zero scalar matrix in $\Mn$.
\end{lem}

\begin{proof} First note that if $A_i$ and $A_j$ are
  in the same $\PGLn$-orbit then $s(A_i) = s(A_j)$. Hence, we may
  remove $A_j$ from the set $\{A_1,\ldots,A_r\}$. After repeating this process
  finitely many times, we may assume that no two of the points $A_1,
  \dots, A_r$ lie in the same $\PGLn$-orbit.
  
  By the above-mentioned theorem of Formanek and Razmyslov, there
  exists a central polynomial $c = c(X_1, \dots, X_N)\in G_{N,n}$ for
  $n \times n$-matrices.  Choose $b_1, \dots, b_N \in \Mn$ such that
  $c(b_1, \dots, b_N) \ne 0$.  We now define $s$ by modifying $c$ as
  follows:
  \[ s(X_1, \dots, X_m) = c\bigl(p_1(X_1, \dots, X_m), \dots, p_N(X_1,
     \dots, X_m)\bigr) \, , \] 
  where the elements $p_j = p_j(X_1, \dots, X_m) \in \Gmn$ will be
  chosen below so that for every $j = 1, \dots, N$,
  \begin{equation} \label{e.p_j}
    p_j(A_1) = p_j(A_2) = \dots = p_j(A_r) = b_j \in \Mn \, .
  \end{equation}
  We first check that this polynomial has the desired properties.
  Being an evaluation of a central polynomial for $n\times n$
  matrices, $s$ is a central element in $\Gmn$ and a polynomial
  identity for all prime $k$-algebras of PI-degree $n-1$.  Moreover,
  $s(A_i) = c(b_1, \dots, b_N) \ne 0$ for every $i = 1, \dots, r$.
  Thus $s$ itself is a central polynomial for $n\times n$ matrices.
  Consequently, $s(A_i)$ is a central element in $\Mn$, i.e., a scalar
  matrix.
  
  It remains to show that $p_1, \dots, p_N \in G_{m, n}$ can be
  chosen so that~\eqref{e.p_j} holds. Consider the representation
  $\phi_i \colon \Gmn \lra \Mn$ given by $p \longmapsto p(A_i)$.  Since each
  $A_i$ lies in $U_{m, n}$, each $\phi_i$ is surjective.  Moreover, by
  our assumption on $A_1, \dots, A_r$, no two of them are conjugate
  under $\PGLn$, i.e., no two of the representations $\phi_i$ are
  equivalent.  The kernels of the $\phi_i$ are thus pairwise distinct
  by \cite[Theorem (9.2)]{artin1}.  Hence the Chinese Remainder Theorem
  tells us that $\phi_1 \oplus \ldots \oplus \phi_r \colon \Gmn \lra
  (\Mn)^r$ is surjective; $p_j$ can now be chosen to be any preimage
  of $(b_j, \dots, b_j) \in (\Mn)^r$.  This completes the proof of
  Lemma~\ref{prel3.1}.
\end{proof}

\section{Definition and first properties of $n$-varieties}
\label{sect.defn}

\begin{defn} \label{def.n-var1}
  (a) An {\it $n$-variety} $X$ is a closed $\PGL_n$-invariant
  subvariety of $U_{m, n}$ for some $m \ge 2$.  In other words, $X =
  \overline{X} \cap U_{m, n}$, where $\overline{X}$ is the Zariski
  closure of $X$ in $\Mnm$.  Note that $X$ is a generically free
  $\PGLn$-variety (in fact, for every $x\in X$, the stabilizer of $x$
  in $\PGLn$ is trivial).

\smallskip
(b) Given a subset $S \subset G_{m, n}$ (or $S \subset T_{m, n}$),
we define its zero locus as
\[ \mathcal{Z}(S) = \{ a = (a_1, \dots, a_m) \in U_{m, n} \,
| \, p(a) = 0, \; \; \forall
p \in S \} \, . \]
Of course, $\calZ(S) = \calZ(J)$, where $J$ is the 2-sided
ideal of $\Gmn$ (or $T_{m, n}$) generated by $S$.
Conversely, given an $n$-variety
$X \subset U_{m, n}$ we define its ideal as
\begin{align*}
  \mathcal{I}(X) &= \{ p \in \Gmn \, | \, p(a) = 0, \; \; \forall
  a \in X \} \, .\\
\intertext{Similarly we define the ideal of $X$ in $\Tmn$, as}
  \mathcal{I}_T(X)&=\{ p \in \Tmn \, | \, p(a) = 0, \; \; \forall
a \in X \} \, .
\end{align*}
Note that $\calI(X)=\calI_T(X) \cap\Gmn$.

\smallskip (c) The {\it polynomial identity coordinate ring} (or {\it
  PI-coordinate ring}) of an $n$-variety $X$ is defined as $G_{m,
  n}/\mathcal{I}(X)$. We denote this ring by $k_n[X]$.
\end{defn}

\begin{remark} \label{rem.regular}
  Elements of $k_n[X]$ may be viewed as
  $\PGLn$-equi\-vari\-ant morphisms $X \lra \Mn$.  The example
  below shows that not every $\PGLn$-equiva\-riant
  morphism $X \lra \Mn(k)$ is of this form.
  On the other hand, if $X$ is irreducible, we will later
  prove that every $\PGLn$-equivariant
  {\em rational} map $X \dasharrow \Mn(k)$
  lies in the total ring of fractions of $k_n[X]$;
  see Proposition~\ref{prop.rat2}.
\end{remark}

\begin{example}\label{ex.regular}
Recall that $U_{2, 2}$ is the open subset of $\M_{2, 2}$
defined by the inequality $c(X_1, X_2) \ne 0$, where
\begin{align*} c(X_1, X_2)  =  \big( 2 \tr(X_1^2)  - \tr(X_1)^2 \big)
\big( 2 \tr(X_2^2)  - \tr(X_2)^2 \big) - \\ \big( 2
 \tr(X_1 X_2) -  \tr(X_1) \det(X_2) \big)^2 \, ; 
\end{align*}
see, e.g.,~\cite[p. 191]{friedland}.
Thus for $X=U_{2,2}$, the $\PGLn$-equivariant morphism
$f \colon X\lra\M_2$ given by
$(X_1, X_2) \longmapsto \frac{1}{c(X_1, X_2)} I_{2 \times 2}$
is not in $k_2[X]=G_{2,2}$ (and not even in $T_{2,2}$).
\end{example}

\begin{remark} \label{rem.points}
  (a) Let $J$ be an ideal of $G_{m, n}$.  Then the points of
  $\mathcal{Z}(J)$ are in bijective correspondence with the surjective
  $k$-algebra homomorphisms $\phi \colon G_{m, n} \lra \Mn$ such that
  $J \subset \Ker(\phi)$ (or equivalently, with the surjective
  $k$-algebra homomorphisms $G_{m,n}/J \lra \Mn$).  Indeed, given $a
  \in \mathcal{Z}(J)$, we associate to it the homomorphism $\phi_a$
  given by $\phi_a \colon p \mapsto p(a)$.
  Conversely, a surjective homomorphism $\phi
  \colon G_{m, n} \lra \Mn$ such that $J\subset\Ker(\phi)$
  gives rise to the point
  \[ a_\phi = (\phi(X_1), \dots, \phi(X_m))
  \in \mathcal{Z}(J) \, , \]
  where $X_i \in G_{m, n}$ is the $i$-th generic matrix in $G_{m, n}$.
  One easily checks that the assignments $a\mapsto\phi_a$ and
  $\phi\mapsto a_\phi$ are inverse to each other.

  (b) The claim in part (a) is also true for an ideal $J$ of $T_{m,
    n}$.  That is, the points of $\mathcal{Z}(J)$ are in bijective
  correspondence with the surjective $k$-algebra homomorphisms $\phi
  \colon T_{m, n} \lra \Mn$ such that $J \subset \Ker(\phi)$ (or
  equivalently, with the surjective $k$-algebra homomorphisms
  $T_{m,n}/J \lra \Mn$).  The proof goes through without changes.
\end{remark}

\begin{lem} \label{lem1.05} Let $a = (a_1, \dots, a_m) \in \Umn$ 
and let $J$ be an ideal of $\Gmn$ {\upshape(}or of
$T_{mn}${\upshape)}. Let 
\[ J(a) = \{ j(a) \, | \, j \in J \} \subset \Mn \, . \] 
Then either $J(a) = (0)$ {\upshape(}i.e, $a \in \calZ (J)${\upshape)}
or  $J(a) = \Mn$.
\end{lem}

\begin{proof} Since $a_1, \dots, a_m$ generate $\Mn$,
$\phi_a(J)$ is a (2-sided) ideal of $\Mn$. Since
$\Mn$ is simple, the lemma follows.
\end{proof}

\begin{lem} \label{lem1.1}
  {\upshape(a)} $\calZ(J) = \calZ(J \cap G_{m, n})$ for
  every ideal $J \subset T_{m, n}$.

  \smallskip {\upshape(b)} If $X \subset U_{m, n}$ is an
  $n$-variety, then $X=\mathcal{Z}(\mathcal{I}(X)) =
  \calZ(\calI_T(X))$.
\end{lem}

\begin{proof} (a) Clearly, $\calZ(J) \subset \calZ(J \cap G_{m, n})$.
  To prove the opposite inclusion, assume the contrary: there exists a
  $y \in U_{m, n}$ such that $p(y) = 0$ for every $p \in J \cap G_{m,
    n}$ but $f(y) \ne 0$ for some $f \in J$. By Lemma~\ref{prel3.1}
  there exists a central polynomial $s \in \Gmn$ for $n\times
  n$-matrices such that $s(y) \ne 0$. By
  \cite[Theorem~1]{schelter:78}, $p = s^if$ lies in $\Gmn$ (and hence,
  in $J \cap G_{m, n}$) for some $i \ge 0$.  Our choice of $y$ now
  implies $0 = p(y) = s^i(y) f(y)$.  Since $s(y)$ is a non-zero
  element of $k$, we conclude that $f(y) = 0$, a contradiction.

  (b) Clearly $X \subset \calZ(\calI_T(X)) \subset \calZ(\calI(X))$.
  Part (a) (with $J = \calI_T(X)$) tells us that $\calZ(\calI(X)) =
  \calZ(\calI_T(X))$.  It thus remains to be shown that
  $\calZ(\calI_T(X)) \subset X$.  Assume the contrary: there exists a
  $z \in U_{m, n}$ such that $p(z) = 0$ for every $p \in\calI_T(X)$
  but $z \not \in X$. Since $X = \overline{X} \cap U_{m, n}$, where
  $\overline{X}$ is the closure of $X$ in $\Mnm$, we conclude that $z
  \not \in \overline{X}$. Let $C = \PGLn \cdot z$ be the orbit of $z$
  in $\Mnm$. Since $z \in U_{m, n}$, $C$ is closed in $\Mnm$;
  see Proposition~\ref{prel.matr-inv}(b). Thus $C$ and
  $\overline{X}$ are disjoint closed $\PGLn$-invariant subsets
  of $\Mnm$. By~\cite[Corollary 1.2]{mf}, there exists a
  $\PGLn$-invariant regular function $f \colon \Mnm \lra k$ such that
  $f \equiv 0$ on $\overline{X}$ but $f \not \equiv 0$ on $C$.  The
  latter condition is equivalent to $f(z) \ne 0$.  Identifying
  elements of $k$ with scalar matrices in $\Mn$, we may view $f$ as a
  central element of $T_{m, n}$.  So $f \in \mathcal{I}_T(X)$ but
  $f(z) \ne 0$, contradicting our assumption.
\end{proof}

\section{Irreducible $n$-varieties}
Of particular interest to us will be {\em irreducible} $n$-varieties.
Here ``irreducible" is understood with respect to the $n$-Zariski
topology on $U_{m, n}$, where the closed subsets are the
$n$-varieties. However, since $\PGL_n$ is a connected group, each
irreducible component of $X$ in the usual Zariski topology is
$\PGLn$-invariant. Consequently, $X$ is irreducible in the $n$-Zariski
topology if and only if it is irreducible in the usual Zariski
topology.

\begin{lem} \label{lem.irred}
Let $\emptyset \neq X \subset U_{m, n}$ be an
$n$-variety. The following are equivalent:
\begin{itemize}
\item[\rm(a)]$X$ is irreducible.
\item[\rm(b)]$\calI_T(X)$ is a prime ideal of $\Tmn$.
\item[\rm(c)]$\calI(X)$ is a prime ideal of $\Gmn$.
\item[\rm(d)]$k_n[X]$ is a prime ring.
\end{itemize}
\end{lem}

\begin{proof}
  (a) $\Rightarrow$ (b): Suppose  that $\mathcal{I}_T(X)$ is not
  prime, i.e., there are ideals $J_1$ and $J_2$ such that $J_1 \cdot J_2
  \subset \mathcal{I}_T(X)$ but $J_1, J_2 \not \subset
  \mathcal{I}_T(X)$.  We claim that $X$ is not irreducible. 
  Indeed, by Lemma~\ref{lem1.05}, $X \subset \calZ(J_1) \cup \calZ(J_2)$.
  In other words, we can write $X = X_1 \cup X_2$,
  as a union of two $n$-varieties,
  where $X_1 = \mathcal{Z}(J_1) \cap X$ and $X_2 =
  \mathcal{Z}(J_2) \cap X$.  It remains to be shown that
  $X_i \ne X$ for $i =
  1, 2$. Indeed, if say, $X_1 = X$ then every element of $J_1$
  vanishes on all of $X$, so that $J_1 \subset \mathcal{I}_T(X)$,
  contradicting our assumption.

 (b) $\Rightarrow$ (c): Clear, since $\calI(X)=\calI_T(X)\cap\Gmn$.

 (c) $\Leftrightarrow$ (d): $\calI(X) \subset G_{m, n}$ is, by
 definition, a prime ideal if and only if $k_n[X]=\Gmn/\calI(X)$ is a
 prime ring.

 (c) $\Rightarrow$ (a): Assume $\calI(X)$ is prime and
  $X = X_1 \cup X_2$ is a union of two $n$-varieties in $\Umn$.
  Our goal is to show that $X = X_1$ or $X = X_2$. Indeed, 
  $\calI(X_1) \cdot \calI(X_2) \subset \calI(X)$ implies
  $\calI(X_i) \subset \calI(X)$ for $i = 1$ or $2$. Taking the zero
  loci and using Lemma~\ref{lem1.1}(b), we obtain
\[ X_i = \calZ(\calI(X_i)) \supset  \calZ(\calI(X)) = X \, , \]
as desired.
\end{proof}

\begin{prop} \label{prop.central.Tmn}
  Let $J \in \Spec_n(T_{m, n})$. Then
  \begin{itemize}
  \item[\rm(a)]$\mathcal{Z}(J) = \mathcal{Z}(J \cap C_{m, n})$
  \item[\rm(b)]$\mathcal{Z}(J)$ is irreducible.
  \end{itemize}
\end{prop}

\begin{proof}
(a) Clearly $\mathcal{Z}(J \cap C_{m, n}) \supset \mathcal{Z}(J)$.
To prove the opposite inclusion, suppose $a \in  
\mathcal{Z}(J \cap C_{m, n})$ and consider the evaluation map
$\phi_a \colon \Tmn \lra \Mn$ given by $\phi_a(p) = p(a)$.

Recall that $\phi_a$ is trace-preserving (see, e.g.,
\cite[Theorem~2.2]{amitsur:small}). Since $\tr(j) \in J \cap \Cmn$
for every $j \in J$, we see that $\tr(j(a)) = 0$ for every $j \in J$.
By Lemma~\ref{lem1.05} this implies that $a \in \calZ(J)$, as claimed.

 \smallskip
 (b) Consider the categorical quotient map $\pi \colon \Mnm \lra Q_{m, n}$
 for the $\PGLn$-action on $\Mnm$. Recall that $C_{m, n} = k[Q_{m, n}]$
 is the coordinate ring of $Q_{m, n}$.
 Note that elements of $C_{m, n}$ may be viewed in two
 ways: as regular functions on $Q_{m, n}$ or (after
 composing with $\pi$) as a $\PGLn$-invariant regular
 function on $\Mnm$. Let $Y \subset Q_{m, n}$
 be the zero locus of $J \cap C_{m, n}$
 in $Q_{m, n}$. Then by part (a),
 \[\mathcal{Z}(J) = \mathcal{Z}(J\cap\Cmn)=\pi^{-1}(Y) \cap U_{m, n} \, .\]
 Since $J$ is a prime ideal of $T_{m, n}$, $J \cap C_{m, n}$ is
 a prime ideal of $C_{m, n}$; see, e.g., \cite[Theorem II.6.5(1)]{procesi2}.
 Hence, $Y$ is irreducible. Now by Proposition~\ref{prel.matr-inv}(e),
 we conclude that $\calZ(J) = \pi^{-1}(Y) \cap U_{m,n}$ is also
 irreducible, as claimed.
\end{proof}

\begin{cor} \label{cor.central.Gmn}
If $J_0 \in \Spec_n(G_{m, n})$, then $\mathcal{Z}(J_0)$ is irreducible.
\end{cor}

\begin{proof} By Lemma~\ref{lem:TraceRingofCoordRing:new},
$J_0 = J \cap \Gmn$ for some $J \in \Spec_n(\Tmn)$. By
Lemma~\ref{lem1.1}(a),
$\calZ(J_0) = \calZ(J)$, and by Proposition~\ref{prop.central.Tmn}(b),
$\calZ(J)$ is irreducible.
\end{proof}

\section{The Nullstellensatz for prime ideals}
\label{sect.nullstellensatz}

\begin{prop} \label{prop.weak-null.prime}
{\upshape (Weak form of the Nullstellensatz)}
  Let $A$ denote the algebra $\Gmn$ or $\Tmn$, and let $J$ be a prime
  ideal of $A$.  Then $\calZ(J)\neq\emptyset$ if and only if $A/J$ has
  PI-degree~$n$.
\end{prop}

Note that for $n = 1$, Proposition~\ref{prop.weak-null.prime} reduces
to the usual (commutative) weak Nullstellensatz (which is used in the
proof of Proposition~\ref{prop.weak-null.prime}). Indeed, a prime ring
of PI-degree~$1$ is simply a nonzero commutative domain; in this case
$G_{m, 1}/J = k[x_1, \dots, x_m]/J$ has PI-degree~1 if and only if $J
\ne k[x_1, \dots, x_m]$.

\begin{proof}
  First assume that $\mathcal{Z}(J) \ne \emptyset$.  Since $A$ has
  PI-degree $n$, its quotient $A/J$ clearly has PI-degree $\le n$. To
  show $\PIdeg(A/J) \ge n$, recall that a point $a = (a_1, \dots, a_m)
  \in Z(J)$ gives rise to a surjective $k$-algebra homomorphism $A/J
  \lra \Mn$; see Remark~\ref{rem.points}.
  
  Conversely, assume that $R=A/J$ is a $k$-algebra of PI-degree $n$.
  Note that $R$ is a Jacobson ring (i.e., the intersection of its
  maximal ideals is zero), and that it is a Hilbert $k$-algebra
  (i.e., every simple homomorphic image is finite-dimensional over $k$
  and thus a matrix algebra over $k$), see
  \cite[Corollary~1.2]{amitsur:procesi}.  So if $c$ is a nonzero
  evaluation in $R$ of a central polynomial for $n\times n$-matrices,
  there is some maximal ideal $M$ of $R$ not containing~$c$.  Then
  $R/M\simeq \Mn$, and we are done in view of Remark~\ref{rem.points}.
\end{proof}
\begin{cor} \label{cor:trace-ring}
  For any irreducible $n$-variety $X$, $\Tmn/\calI_T(X)$ is the trace
  ring of the prime $k$-algebra $k_n[X]$.
\end{cor}

\begin{proof}
  By Lemma~\ref{lem.irred} and Proposition~\ref{prop.weak-null.prime},
  $\calI_T(X)$ is a prime ideal of $\Tmn$ of PI-degree~$n$.
  Consequently, $\calI(X) = \calI_T(X) \cap \Gmn$ is a prime ideal of
  $\Gmn$, and the desired conclusion follows from
  Lemma~\ref{lem:TraceRingofCoordRing:new}(a).
\end{proof}

\begin{prop} \label{prop.strong-null}
{\upshape (Strong form of the Nullstellensatz)}
\begin{itemize}
\item[\upshape(a)]$\mathcal{I}(\mathcal{Z}(J_0)) = J_0$ for every $J_0
  \in \Spec_n(G_{m, n})$.
\item[\upshape(b)] $\mathcal{I}_T(\mathcal{Z}(J)) = J$ for every $J
  \in \Spec_n(T_{m, n})$.
\end{itemize}
\end{prop}

For $n = 1$ both parts reduce to the usual (commutative)
strong form of the Nullstellensatz for prime ideals (which is used in
the proof of Proposition~\ref{prop.strong-null}).

\begin{proof}
We begin by reducing part (a) to part (b).
Indeed, by Lemma~\ref{lem:TraceRingofCoordRing:new},
$J_0 = J \cap G_{m, n}$ for some $J \in \Spec_n(T_{m, n})$.
Now
\[\calI(\calZ(J_0))
\stackrel{\text{\tiny (1)}}{=}
\calI(\calZ(J))=\calI_T(\calZ(J))\cap\Gmn
\stackrel{\text{\tiny (2)}}{=} J \cap\Gmn=J_0 \,,\]
where (1) follows from  Lemma~\ref{lem1.1}(a) and (2) follows from
part~(b).

It thus remains to prove (b).  Let $X = \mathcal{Z}(J)$.
Then $X \ne \emptyset$ (see Proposition~\ref{prop.weak-null.prime}),
$X$ is irreducible (see Proposition~\ref{prop.central.Tmn}(b)) and
hence $\mathcal{I}_T(X)$ is a prime ideal of $T_{m, n}$
(see Lemma~\ref{lem.irred}).  Clearly $J \subset \mathcal{I}_T(X)$;
our goal is to show that $J = \mathcal{I}_T(X)$. In fact, we only
need to check that
\begin{equation} \label{e.strong-null}
J \cap C_{m, n} = \mathcal{I}_T(X) \cap C_{m, n} \, .
\end{equation}
Indeed, suppose~\eqref{e.strong-null} is established.
Choose $p \in \mathcal{I}_T(X)$; we want to show that $p \in J$.
For every $q \in T_{m, n}$ we have $pq \in I_T(X)$ and thus
\[ \tr(p \cdot q) \in \mathcal{I}_T(X) \cap C_{m, n} = J \cap C_{m, n}
\, , \]
see Lemma~\ref{lem:TraceRingofCoordRing:new}(b).
Hence, if we denote the images of $p$ and $q$ in $T_{m, n}/J$
by $\overline{p}$ and $\overline{q}$ respectively,
Lemma~\ref{lem:TraceRingofCoordRing:new}(a) tells us that
$\tr(\overline{p} \cdot \overline{q}) = 0$ in $T_{m, n}/J$ for
every $\overline{q} \in T_{m,n}/J$.
Consequently, $\overline{p} = 0$, i.e., $p \in J$, as desired.

We now turn to proving~\eqref{e.strong-null}.
Consider the categorical quotient map $\pi \colon \Mnm \lra Q_{m, n}$
for the $\PGLn$-action on $\Mnm$; here $C_{m, n} = k[Q_{m, n}]$
is the coordinate ring of $Q_{m, n}$.
Given an ideal $H \subset C_{m, n}$, denote its zero locus
in $Q_{m, n}$ by
\[ \mathcal{Z}_0(H) = \{ a \in Q_{m, n} \, | \, h(a) = 0 \; \forall \,
h \in H \} \, . \]
Since $\mathcal{Z}(J \cap C_{m, n}) = \mathcal{Z}(J) = X$ in $U_{m, n}$
(see Proposition~\ref{prop.central.Tmn}(a)), we conclude that
$\mathcal{Z}_0(J \cap C_{m, n}) \cap \pi(U_{m, n}) = \pi(X)$.
On the other hand, since $\pi(U_{m, n})$ is Zariski open in
$Q_{m, n}$ (see Lemma~\ref{prel.matr-inv}(d)),
and $J \cap C_{m, n}$ is a prime ideal of $C_{m, n}$
(see~\cite[Theorem II.6.5(1)]{procesi2}), we have
\begin{equation} \label{e1.strong-null}
\mathcal{Z}_0(J \cap C_{m, n}) = \overline{\pi(X)} \, ,
\end{equation}
where $\overline{\pi(X)}$ is the Zariski closure of $\pi(X)$
in $Q_{m, n}$.
Now suppose $f \in \mathcal{I}_T(X) \cap C_{m, n}$. Our goal is to show that
$f \in J \cap C_{m, n}$.  Viewing $f$ as an element of $C_{m, n}$,
i.e., a regular function on $Q_{m, n}$,
we see that $f \equiv 0$ on $\pi(X)$ and hence, on $\overline{\pi(X)}$.
Now applying the usual (commutative) Nullstellensatz to the prime ideal
$J \cap C_{m, n}$ of $C_{m, n}$, we see that~\eqref{e1.strong-null}
implies $f \in J \cap C_{m, n}$, as desired.
\end{proof}

\begin{remark}\label{rem.Amitsur}
  In this paper, we consider zeros of ideals of $\Gmn$ in $\Umn$.  In
  contrast, Amitsur's Nullstellensatz~\cite{amitsur} (see
  also~\cite{amitsur:procesi}) deals with zeros in the larger space
  $\Mnm$.  Given an ideal $J$ of $\Gmn$, denote by $\calZ(J;\Mnm)$ the
  set of zeroes of $J$ in $\Mnm$.  Since
  $\calZ(J;\Mnm)\supset\calZ(J)$, it easily follows that
  \[ J \subset \calI\bigl(\calZ(J;\Mnm)\bigr)
       \subset\calI\bigl(\calZ(J)\bigr)\,.
  \]
  One particular consequence of Amitsur's Nullstellensatz is that the
  first inclusion is an equality if $J$ is a prime ideal.
  Proposition~\ref{prop.strong-null}  implies that both inclusions are
  equalities, provided $J$ is a prime ideal of PI-degree~$n$.  Note
  that the second inclusion can be strict, e.g., if $J$ is a prime
  ideal of PI-degree $<n$ (since then $\calZ(J)=\emptyset$ by
  Proposition~\ref{prop.weak-null.prime}). 
\end{remark}

The following theorem summarizes many of our results so far.

\begin{thm}\label{thm.inverse.bijections}
Let $n \ge 1$ and $m \ge 2$ be integers.
\begin{itemize}
\item[\upshape(a)]$\calZ(\hbox to.6em{\hrulefill})$ and $\calI(\hbox
  to.6em{\hrulefill})$ are mutually inverse inclusion-reversing
  bijections between $\Spec_n(\Gmn)$ and the set of irreducible
  $n$-varieties $X \subset U_{m, n}$.
\item[\upshape(b)]$\calZ(\hbox to.6em{\hrulefill})$ and $\calI_T(\hbox
  to.6em{\hrulefill})$ are mutually inverse inclusion-reversing
  bijections between $\Spec_n(\Tmn)$ and the set of irreducible
  $n$-varieties $X \subset U_{m, n}$.  \qed
\end{itemize}
\end{thm}

\section{Regular maps of $n$-varieties}
\label{sect.reg}

Recall that an element $g$ of $G_{m, n}$ may be viewed as a regular
$\PGLn$-equi\-variant map $g \colon \Mnm \lra \Mn$.  Now suppose $X$
is an $n$-variety in $\Umn$. Then, restricting $g$ to $X$, we see that
$g|_X = g'|_X$ if and only if $g' - g \in \calI(X)$. Hence, elements
of the PI-coordinate ring of $X$ may be viewed as $\PGLn$-equivariant
morphisms $X \lra \Mn$. All of this is completely analogous to the
commutative case, where $n = 1$, $\M_1 = k$, $k_1[X] = k[X]$ is the
usual coordinate ring of $X \subset k^m$, and elements of $k[X]$ are
the regular functions on $X$. It is thus natural to think of elements
of $k_n[X]$ as ``regular functions" on $X$, even though we shall not
use this terminology.  (In algebraic geometry, functions are usually
assumed to take values in the base field~$k$, while elements of
$k_n[X]$ take values in $\Mn$.)
We also remark that not every $\PGLn$-equivariant morphism
$X \lra \Mn$ of algebraic varieties (in the usual sense)
is induced by elements of $k_n[X]$, see Example~\ref{ex.regular}.

\begin{defn} \label{def.n-var2}
Let $X \subset U_{m, n}$ and $Y \subset U_{l, n}$ be $n$-varieties.

\smallskip (a) A map $f \colon X \lra Y$ is called a
{\it regular map} of $n$-varieties, if it is of the form
$f=(f_1,\ldots,f_l)$ with each $f_i\in k_n[X]$ (so that $f$ sends $a =
(a_1, \dots, a_{m}) \in X \subset \Mnm$ to $(f_1(a), \dots, f_l(a))
\in Y \subset \Mnl$).  Note that a morphism of $n$-varieties $X\lra Y$
extends to a $\PGLn$-equivariant morphism $\Mnm\lra(\Mn)^l$.

\smallskip
(b) The $n$-varieties $X$ and $Y$ are called {\em isomorphic} if
there are mutually inverse regular maps $X \lra Y$ and $Y \lra X$ .

\smallskip (c) A regular map $f = (f_1, \dots, f_l) \colon X \lra Y
\subset U_{l, n}$ of $n$-varieties induces a $k$-algebra homomorphism
$f^*\colon k_n[Y] \lra k_n[X]$ given by $\overline{X_i} \lra f_i$ for
$i = 1, \dots, l$, where $\overline{X_1}, \dots, \overline{X_l}$ are
the images of the generic matrices $X_1, \dots, X_l \in G_{l, n}$ in
$k_n[Y] = G_{l, n}/\mathcal{I}(Y)$.  One easily verifies that for
every $g\in k_n[Y]$, $f^*(g)=g\circ f \colon X \lra Y \lra \Mn$.

\smallskip (d) Conversely, a $k$-algebra homomorphism $\alpha \colon
k_n[Y] \lra k_n[X]$ induces a regular map $\alpha_* = (f_1, \dots,
f_l) \colon X \lra Y$ of $n$-varieties, where $f_i =
\alpha(\overline{X_i})$.  It is easy to check that for
every $g \in k_n[Y]$, $\alpha(g) = g \circ \alpha_* \colon X \lra Y \lra \Mn$.
\end{defn}

\begin{remark} \label{rem.**}
  It is immediate from these definitions that $(f^*)_* = f$ for any
  regular map $f \colon X \lra Y$ of $n$-varieties, and
  $(\alpha_*)^* = \alpha$ for any $k$-algebra homomorphism $\alpha
  \colon k_n[Y] \lra k_n[X]$. Note also that
  $(\id_X)^*=\id_{k_n[X]}$, and $(\id_{k_n[X]})_*=\id_X$.
\end{remark}

\begin{lem} \label{lem.reg1}
  Let $X \subset U_{m, n}$ and $Y \subset U_{l, n}$ be $n$-varieties,
  and let $k_n[X]$ and $k_n[Y]$ be their respective PI-coordinate
  rings.
  \begin{itemize}
  \item[\rm(a)]If $f \colon X \lra Y$ and $g \colon Y \lra Z$ are
    regular maps of $n$-varieties, then $(g \circ f)^* = f^* \circ g^*$.
  \item[\rm(b)]If $\alpha \colon k_n[Y] \lra k_n[X]$ and $\beta \colon
    k_n[Z] \lra k_n[Y]$ are $k$-algebra homomorphisms, then $(\alpha
    \circ \beta)_* = \beta_* \circ \alpha_*$.
  \item[\rm(c)]$X$ and $Y$ are isomorphic as $n$-varieties if and only
    if $k_n[X]$ and $k_n[Y]$ are isomorphic as $k$-algebras.
  \end{itemize}
\end{lem}

\begin{proof} Parts (a) and (b) follow directly from
  Definition~\ref{def.n-var2}. The proofs are exactly the same as in
  the commutative case (where $n = 1$); we leave them as an exercise
  for the reader.

To prove (c), suppose
$f \colon X \lra Y$ and $g \colon Y \lra X$
are mutually inverse morphisms of $n$-varieties.
Then by part (a), $f^* \colon k_n[Y] \lra k_n[X]$
and $g^* \colon k_n[X] \lra k_n[Y]$ are mutually inverse
$k$-algebra homomorphisms, showing
that $k_n[X]$ and $k_n[Y]$ are isomorphic.

Conversely, if $\alpha \colon k_n[Y] \lra k_n[X]$ and
$\beta \colon k_n[X] \lra k_n[Y]$ are mutually inverse
homomorphisms of $k$-algebras then by part (b),
$\alpha_*$ and $\beta_*$ are mutually inverse morphisms
between the $n$-varieties $X$ and $Y$.
\end{proof}

\begin{thm} \label{thm.n-Var}
  Let $R$ be a finitely generated prime $k$-algebra of PI-degree~$n$.
  Then $R$ is isomorphic {\upshape(}as a $k$-algebra{\upshape)} to
  $k_n[X]$ for some irreducible $n$-variety $X$.  Moreover, $X$ is
  uniquely determined by $R$, up to isomorphism of $n$-varieties.
\end{thm}

\begin{proof}
  By our assumptions on $R$ there exists a surjective ring
  homomorphism $\varphi \colon G_{m, n} \lra R$.  Then $J_0 =
  \Ker(\varphi)$ lies in $\Spec_n(G_{m, n})$.  Set $X =
  \mathcal{Z}(J_0) \subset U_{m, n}$. Then $X$ is irreducible (see
  Corollary~\ref{cor.central.Gmn}), and $J_0 = \mathcal{I}(X)$ (see
  Proposition~\ref{prop.strong-null}(a)).  Hence $R$ is isomorphic to
  $G_{m, n}/\mathcal{I}(X) = k_n[X]$, as claimed.  The uniqueness of
  $X$ follows from Lemma~\ref{lem.reg1}(c).
\end{proof}

We are now ready to prove Theorem~\ref{thm1}. Recall that
for $n = 1$, Theorem~\ref{thm1} reduces to~\cite[Corollary 3.8]{hartshorne}.
Both are proved by the same argument.  Since the proof
of~\cite[Corollary 3.8]{hartshorne} is omitted in~\cite{hartshorne},
we reproduce this argument here for the sake of completeness.

\smallskip
\begin{proof}[Proof of Theorem~\ref{thm1}]
  By Lemma~\ref{lem.reg1}, the contravariant functor~$\mathcal F$ in
  Theorem~\ref{thm1} is well-defined.  It is full and faithful by
  Remark~\ref{rem.**}.  Moreover, by Theorem~\ref{thm.n-Var}, every
  object in $\PI_n$ is isomorphic to the image of an object in
  $\Var_n$.  Hence $\mathcal F$ is a covariant equivalence of
  categories between $\Var_n$ and the dual category of $\PI_n$, see,
  e.g., \cite[Theorem~7.6]{blyth}.  In other words, $\mathcal F$ is
  a contravariant equivalence of categories between
  $\Var_n$ and $\PI_n$, as claimed.
\end{proof}

We conclude our discussion of regular maps of $n$-varieties with an
observation which we will need in the next section.

\begin{lem}\label{lem:central-in-kn[X]}
  Let $X$ be an irreducible $n$-variety, and $c$ a central element of
  $k_n[X]$ or of its trace ring $\Tmn/\calI_T(X)$.  Then the image of $c$
  in $\Mn$ consists of scalar matrices.
\end{lem}

\begin{proof}
  Denote by $k_n(X)$ the common total ring of fractions of $k_n[X]$
  and $\Tmn/\calI_T(X)$.  By Lemma~\ref{prel3.1}, there exists a central
  polynomial $s\in\Gmn$ for $n\times n$ matrices which
  does not identically vanish on $X$.  Then $R=\Gmn[s\inv]$ is an
  Azumaya algebra.  Consider the natural map $\phi\colon R \lra
  k_n(X)$.  Then the center of $\phi(R)$ is $\phi(\text{Center}(R))$;
  see, e.g., \cite[Proposition 1.11]{demeyer}.
  Note that $k_n(X)$ is a central localization of $\phi(R)$.  Hence
  the central element $c$ of $k_n(X)$ is of the form
  $c=\phi(p)\phi(q)\inv$ for central elements $p,q\in\Gmn$ with
  $q\not\equiv0$ on~$X$.  All images of $p$ and $q$ in $\Mn$ are
  central, i.e., scalar matrices.  Thus $c(x)$ is a scalar matrix for
  each $x$ in the dense open subset of~$X$ on which $q$ is nonzero.
  Consequently, $c(x)$ is a scalar matrix for every $x\in X$.
\end{proof}

\section{Rational maps of $n$-varieties}
\label{sect.ratl}

\begin{defn}\label{def.ratl.1}
  Let $X$ be an irreducible $n$-variety.
  The total ring of fractions of the prime algebra $k_n[X]$ will
  be called the {\it central simple algebra of rational
  functions on $X$} and denoted by $k_n(X)$.
\end{defn}

\begin{remark}\label{rem.rational.Tmn}
  One can also define $k_n(X)$ using the trace ring instead
  of the generic matrix ring. That is, $k_n(X)$ is also
  the total ring of fractions of $\Tmn/\calI_T(X)$. Indeed, by
  Corollary~\ref{cor:trace-ring}, $\Tmn/\calI_T(X)$ is the
  trace ring of $k_n[X]$, so the two have the same total
  ring of fractions.
\end{remark}

Recall that $k_n(X)$ is obtained from $k_n[X]$ by inverting all
non-zero central elements; see, e.g.,~\cite[Theorem 1.7.9]{rowen:PI}.
In other words, every $f \in k_n(X)$ can be written as $f = c^{-1} p$,
where $p \in k_n[X]$ and $c$ is a nonzero central element of $k_n[X]$.
Recall from Lemma~\ref{lem:central-in-kn[X]} that for each $x\in X$,
$c(x)$ is a scalar matrix in $\Mn$, and thus invertible if it is
nonzero.  Viewing $p$ and $c$ as $\PGLn$-equivariant morphisms $X \lra
\Mn$ (in the usual sense of commutative algebraic geometry), we see
that $f$ can be identified with a rational map $c^{-1} p \colon X
\dasharrow \Mn$.  One easily checks that this map is independent of
the choice of $c$ and $p$, i.e., remains the same if we replace $c$
and $p$ by $d$ and $q$, such that $f = c^{-1} p = d^{-1} q$. We will
now see that every $\PGLn$-equivariant rational map $X \dasharrow \Mn$
is of this form.

\begin{prop} \label{prop.rat2}
Let $X \subset U_{m, n}$ be an irreducible $n$-variety.
Then the natural inclusion $k_n(X) \hookrightarrow
\RMaps_{\PGLn}(X, \Mn)$ is an isomorphism.
\end{prop}

Here $\RMaps_{\PGLn}(X, \Mn)$ denotes the $k$-algebra of
$\PGLn$-equivariant rational maps $X \dasharrow \Mn$,
with addition and multiplication induced from $\Mn$.
Recall that a regular analogue of Proposition~\ref{prop.rat2}
(with rational maps replaced by regular maps, and $k_n(X)$
replaced by $k_n[X]$) is false; see Remark~\ref{rem.regular}
and Example~\ref{ex.regular}.

\begin{proof}[First proof {\rm(}algebraic{\rm)}]
  Recall that $k_n(X)$ is, by definition, a central simple algebra
  of PI-degree $n$. By~\cite[Lemma 8.5]{reichstein4} (see also
  \cite[Definition 7.3 and Lemma 9.1]{reichstein4}),
  $\RMaps_{\PGLn}(X, \Mn)$ is a central simple algebra
  of PI-degree $n$ as well.  It is thus enough to show that
  the centers of $k_n(X)$ and $\RMaps_{\PGLn}(X, \Mn)$ coincide.

  Let $\overline{X}$ be the closure of $X$ in $\Mnm$.
  By \cite[Lemma 8.5]{reichstein4}, the center of
  $\RMaps_{\PGLn}(X, \Mn) =
  \RMaps_{\PGLn}(\overline{X}, \Mn)$ is the field $k(X)^\PGLn$ of
  $\PGLn$-invariant rational functions $f \colon X \dasharrow k$ (or
  equivalently, the field $k(\overline{X})^{\PGLn}$). Here, as usual,
  we identify $f$ with $f \cdot I_n \colon X \dasharrow \Mn$.
  It now suffices to show that
  \begin{equation} \label{e.rat2.new}
    k(X)^\PGLn = \text{Center}(k_n(X))  \, .
  \end{equation}
  Recall from Lemma~\ref{lem:TraceRingofCoordRing:new}(a) that the natural
  algebra homomorphism $\Gmn\lra k_n(X)$ extends to a homomorphism
  $\Tmn\lra k_n(X)$.  So the center of $k_n(X)$ contains all functions
  $f|_{ X} \colon X \lra k$, as $f$ ranges over the ring $C_{m, n} =
  k[\Mnm]^{\PGLn}$. Since $X \subset U_{m, n}$, these functions
  separate the $\PGLn$-orbits in $X$; see Proposition~\ref{prel.matr-inv}(b).
  Equality~\eqref{e.rat2.new} now follows a theorem 
  of Rosenlicht; cf.~\cite[Lemma 2.1]{pv}.
\end{proof}

\begin{proof}[Alternative proof {\rm (}geometric{\rm )}]
  By Remark~\ref{rem.rational.Tmn}, it suffices to show that
  for every $\PGLn$-equivariant rational map
  $f \colon X \dasharrow \Mn$ there exists an $h \in
  k[\Mn]^{\PGLn}=\Cmn$ such that $h\not\equiv0$ on $X$ and
  $h f \colon x \mapsto h(x) f(x)$ lifts to a regular map $\Mnm \lra
  \Mn$ (and in particular, $hf$ is a regular map $X \lra \Mn$).

  It is enough to show that the ideal $I \subset k[\Mnm]$ given by
  \[ \text{$I = \{ h \in k[(\Mn)^m] \, | \, h f$ lifts to a regular
  map $\Mnm \lra \Mn \}$} \]
  contains a $\PGLn$-invariant element~$h$ such that $h\not\equiv0$ on
  $X$. Indeed, regular $\PGLn$-equivariant morphisms $\Mnm \lra \Mn$ are
  precisely elements of $T_{m, n}$; hence, $hf \colon X \lra \Mn$
  would then lie in $k_n(X)$, and so would $f = h^{\inv} (hf)$, thus
  proving the lemma.

  Denote by $Z$ the zero locus of $I$ in $\Mnm$ (in the usual sense,
  not in the sense of Definition~\ref{def.n-var1}(b)).  Then $Z \cap
  X$ is, by definition the indeterminacy locus of $f$; in particular,
  $X \not \subset Z$.  Choose $a \in X \setminus Z$ and let $C = \PGLn
  \cdot a$ be the orbit of $a$ in $X$. Since $a \in X \subset \Umn$,
  Proposition~\ref{prel.matr-inv}(b) tells us that $C$ is closed
  in $\Mnm$. In summary, $C$ and $Z$ are disjoint $\PGLn$-invariant
  Zariski closed subsets of $\Mnm$. Since $\PGLn$ is
  reductive, they can be separated by a regular invariant, i.e.,
  there exists a $0 \ne j \in k[\Mnm]^{\PGLn}$ such that
  $j(a) \ne 0$ but $j \equiv 0$ on $Z$; see, e.g., \cite[Corollary
  1.2]{mf}.  By Hilbert's Nullstellensatz, $h = j^r$ lies in $I$ for
  some $r \ge 1$.  This $h$ has the desired properties: it is a
  $\PGLn$-invariant element of $I$ which is not identically zero on $X$.
\end{proof}

\begin{defn} \label{def.ratl.2}
  Let $X \subset U_{m, n}$ and $Y \subset U_{l, n}$ be irreducible
  $n$-varieties.

  \smallskip (a) A rational map $f \colon X \dasharrow Y$ is
  called a {\em rational map of $n$-varieties} if
  $f = (f_1,\ldots,f_l)$ where each $f_i\in k_n(X)$.
  Equivalently (in view of Proposition~\ref{prop.rat2}),
  a rational map $X \dasharrow Y$ of $n$-varieties
  is simply a $\PGLn$-equivariant rational map (in the usual sense).

  \smallskip (b) The $n$-varieties $X$ and $Y$ are called {\em
  birationally isomorphic} or {\em birationally equivalent} if there
  exist dominant rational maps of $n$-varieties $f \colon X \dasharrow
  Y$ and $g \colon Y \dasharrow X$ such that $f \circ g = \id_Y$ and
  $g \circ f = \id_X$ (as rational maps of varieties).

  \smallskip (c) A dominant rational map $f=(f_1,\ldots,f_l) \colon X
  \dasharrow Y$ of $n$-varieties induces a $k$-algebra homomorphism
  (i.e., an embedding) $f^* \colon k_n(Y) \lra k_n(X)$ of
  central simple algebras defined by $f^*(\overline{X_i})=f_i$, where
  $\overline{X_i}$ is the image of the generic matrix $X_i\in G_{l,n}$
  in $k_n[Y] \subset k_n(Y)$.  One easily verifies that for every
  $g\in k_n(Y)$, $f^*(g)=g\circ f$, if one views $g$
  as a $\PGLn$-equivariant rational map $Y\dasharrow\Mn$.

  \smallskip (d) Conversely, a $k$-algebra homomorphism (necessarily
  an embedding) of central simple algebras $\alpha \colon k_n(Y) \lra
  k_n(X)$ (over $k$) induces a dominant rational map $f = \alpha_*
  \colon X \dasharrow Y$ of $n$-varieties.  This map is given by $f =
  (f_1, \dots, f_l)$ with $f_i = \alpha(\overline{X_i}) \in k_n(X)$,
  where $\overline{X_1}, \dots, \overline{X_l}$ are the images of the
  generic matrices $X_1, \dots, X_l \in G_{l, n}$.  It is easy to
  check that for every $g\in k_n(Y)$, $\alpha(g)=g\circ\alpha_*$, if
  one views $g$ as a $\PGLn$-equivariant rational map $Y\dasharrow\Mn$.
\end{defn}

\begin{remark} \label{rem.rat5}
  Once again, the identities $(f^*)_* = f$ and $(\alpha_*)^* = \alpha$
  follow directly from these definitions.  Similarly,
$(\id_X)^* = \id_{k_n(X)}$ and
$(\id_{k_n(X)})_* = \id_{X}$.
\end{remark}

We also have the following analogue of
Lemma~\ref{lem.reg1} for dominant rational maps.
The proofs are again the same as in
the commutative case (where $n = 1$);
we leave them as an exercise for the reader.

\begin{lem} \label{lem.rat1}
  Let $X \subset U_{m_1, n}$, $Y \subset U_{m_2, n}$ and $Z \subset
  U_{m_3, n}$ be irreducible $n$-varieties.
  \begin{itemize}
  \item[\rm(a)]If $f \colon X \dasharrow Y$ and $g \colon Y \dasharrow
    Z$ are dominant rational maps of $n$-varieties then $(g \circ f)^*
    = f^* \circ g^*$.
  \item[\rm(b)]If $\alpha \colon k_n(Y) \hookrightarrow k_n(X)$ and
    $\beta \colon k_n(Z) \hookrightarrow k_n(Y)$ are homomorphisms
    {\upshape(}i.e., embeddings{\upshape)} of central simple algebras
    then $(\alpha \circ \beta)_* = \beta_* \circ \alpha_*$.
  \item[\rm(c)]$X$ and $Y$ are birationally isomorphic as
    $n$-varieties if and only if the central simple algebras $k_n(X)$
    and $k_n(Y)$ are isomorphic as $k$-algebras.  \qed
  \end{itemize}
\end{lem}

We are now ready to prove the following birational analogue of
Theorem~\ref{thm.n-Var}.

\begin{thm} \label{thm.W(A)}
  Let $K/k$ be a finitely generated field extension and $A$ be a
  central simple algebra of degree $n$ with center $K$.  Then $A$ is
  isomorphic {\upshape(}as a $k$-algebra{\upshape)} to $k_n(X)$ for some
  irreducible $n$-variety $X$.  Moreover, $X$ is uniquely determined
  by $A$, up to birational isomorphism of $n$-varieties.
\end{thm}

\begin{proof} Choose generators $a_1, \dots, a_N \in K$ for the field
extension $K/k$ and a $K$-vector space basis $b_1, \dots, b_{n^2}$ for $A$.
Let $R$ be the $k$-subalgebra of $A$ generated by all $a_i$ and $b_j$.
By our construction $R$ is a prime $k$-algebra of PI-degree $n$,
with total ring of fraction $A$. By Theorem~\ref{thm.n-Var}
there exists an $n$-variety $X$ such that $k_n[X] \simeq R$
and hence, $k_n(X) \simeq A$. This proves the existence of $X$.
Uniqueness follows from Lemma~\ref{lem.rat1}(c).
\end{proof}

\section{Generically free $\PGLn$-varieties}
\label{sect.gener}

An irreducible $n$-variety is, by definition, an irreducible
generically free $\PGLn$-variety.  The following lemma says that up to
birational isomorphism, the converse is true as well.

\begin{lem} \label{lem.rat3}
  Every irreducible generically free $\PGLn$-variety $X$ is
  birationally isomorphic {\rm(}as $\PGLn$-variety{\rm)} to an
  irreducible $n$-variety in $U_{m, n}$ for some $m \ge 2$.
\end{lem}

\begin{proof} Choose $a \in U_{2, n}$. By \cite[Proposition 7.1]{reichstein4}
  there exists a $\PGLn$-equivariant rational map $\phi \colon X
  \dasharrow (\Mn)^2$ whose image contains $a$. Now choose
  $\PGLn$-invariant rational functions $c_1, \dots, c_{r} \in
  k(X)^{\PGLn}$ on $X$ which separate $\PGLn$-orbits in general
  position (this can be done by a theorem of Rosenlicht; cf.
  e.g.,~\cite[Theorem 2.3]{pv}).  We now set $m = r+2$ and define $f
  \colon X \dasharrow \Mnm$ by
  \[ f(x) = \bigl(c_1(x)I_{n \times n} , \dots, c_{r}(x)I_{n \times
  n}, \phi(x)\bigr) \in (\Mn)^{r} \times (\Mn)^2 = (\M_{n})^m \, . \]
  Let $\overline{Y}$ be the Zariski closure of $f(X)$ in $(\M_n)^m$,
  and $Y = \overline{Y} \cap U_{m, n}$.  By our choice of $\phi$, $Y
  \ne \emptyset$. It thus remains to be shown that $f$ is a birational
  isomorphism between $X$ and $Y$ (or, equivalently, $\overline{Y}$).
  Since we are working over a base field $k$ of characteristic zero,
  it is enough to show that $X$ has a dense open subset $S$ such that
  $f(a) \ne f(b)$ for every pair of distinct $k$-points $a, b \in S$.

  Indeed, choose $S \subset X$ so that (i) the generators $c_1, \dots,
  c_{r}$ of $k(X)^{\PGLn}$ separate $\PGLn$-orbits in $S$, (ii) $f$ is
  well-defined in $S$ and (iii) $f(S) \subset U_{m, n}$. Now let $a$,
  $b \in S$, and assume that $f(a)=f(b)$.  Then $a$ and $b$ must
  belong to the same $\PGLn$-orbit.  Say $b =h(a)$, for some
  $h\in\PGLn$.  Then $f(a)=f(b)=hf(a)h\inv$.  Since $f(a)\in U_{m,
    n}$, $h=1$, so that $a=b$, as claimed.
\end{proof}

Lemma~\ref{lem.rat3} suggests that in the birational setting the
natural objects to consider are arbitrary generically free
$\PGLn$-varieties, rather than $n$-varieties.  The relationship
between the two is analogous to the relationship between affine
varieties and more general algebraic (say, quasi-projective) varieties
in the usual setting of (commutative) algebraic geometry.  In
particular, in general one cannot assign a PI-coordinate ring $k_n[X]$
to an irreducible generically free $\PGLn$-variety in a meaningful
way.  On the other hand, we can extend the definition of
$k_n(X)$ to this setting as follows.

\begin{defn} \label{def.2.kn(X)}
  Let $X$ be an irreducible generically free $\PGLn$-variety.  Then
  $k_n(X)$ is the $k$-algebra of $\PGLn$-equivariant rational maps $f
  \colon X \dasharrow \Mn$, with addition and multiplication induced
  from $\Mn$.
\end{defn}

Proposition~\ref{prop.rat2} tells us that if $X$ is an irreducible
$n$-variety then this definition is consistent with
Definition~\ref{def.ratl.1}.  In place of $k_n(X)$ we will sometimes
write $\RMaps_{\PGLn}(X, \Mn)$.

\begin{defn} \label{def.gener}
  A dominant rational map $f \colon X \dasharrow Y$ of generically
  free $\PGLn$-varieties gives rise to a homomorphism (embedding) $f^*
  \colon k_n(Y) \lra k_n(X)$ given by $f^*(g) = g \circ f$ for every
  $g \in k_n(Y)$.
\end{defn}

If $X$ and $Y$ are $n$-varieties, this definition of $f^*$ coincides
with Definition~\ref{def.ratl.2}(c).  Note that
$(\id_X)^*=\id_{k_n(X)}$, and that $(g\circ f)^*=f^*\circ g^*$ if $g
\colon Y \dasharrow Z$ is another $\PGLn$-equivariant dominant
rational map.  We will now show that Definition~\ref{def.ratl.2} and
Remark~\ref{rem.rat5} extend to this setting as well.

\begin{prop}\label{prop.functor}
  Let $X$ and $Y$ be generically free irreducible $\PGLn$-va\-ri\-eties
  and
  \[\alpha \colon k_n(X) \to k_n(Y)\]
  be a $k$-algebra homomorphism.  Then there is a unique
  $\PGLn$-equivariant, dominant rational map $\alpha_* \colon Y\dasharrow X$
  such that $(\alpha_*)^* = \alpha$.
\end{prop}

\begin{proof} If $X$ and $Y$ are $n$-varieties,
i.e., closed $\PGLn$-invariant subvarieties
of $U_{m, n}$ and $U_{l, n}$ respectively (for some $m, l \ge 2$)
then $\alpha_* \colon Y \dasharrow X$
is given by Definition~\ref{def.ratl.2}(d), and uniqueness
follows from Remark~\ref{rem.rat5}.

In general, Lemma~\ref{lem.rat3} tells us that there are
birational isomorphisms $X \dasharrow X'$ and $Y \dasharrow Y'$
where $X'$ and $Y'$ are $n$-varieties. The proposition is now
a consequence of the following lemma.
\end{proof}

\begin{lem} \label{lem.prop.functor}
Let $f \colon X\dasharrow X'$ and
$g \colon Y\dasharrow Y'$ be birational
isomorphisms of $\PGLn$-varieties. If
Proposition~\ref{prop.functor}
holds for $X'$ and $Y'$ then it holds for $X$ and $Y$.
\end{lem}

\begin{proof}
  Note that by our assumption, the algebra homomorphism
  \[\beta = (g^*)\inv \circ \alpha \circ
  f^* \colon k_n(X') \to k_n(Y')\] is induced by the
  $\PGLn$-equivariant, dominant rational map $\beta_* \colon Y' \to X'$:
  \[\xymatrix{%
    k_n(X) \ar@{->}[r]^{\alpha}
      & k_n(Y) &
      & Y \ar@{-->}[d]_g & X \ar@{-->}[d]^f \\
    k_n(X') \ar@{->}[r]^{\beta} \ar@{->}[u]^{f^*}
      & k_n(Y') \ar@{->}[u]_{g^*} &
      & Y' \ar@{-->}[r]^{\beta_*} & X'
  }\]
Now one easily checks that the dominant rational map
\[\alpha_*:= f\inv \circ \beta_* \circ g \colon Y \dasharrow X\]
has the desired property: $(\alpha_*)^*=\alpha$.
This shows that $\alpha_*$ exists.
To prove uniqueness, let $h \colon Y \lra X$ be
another $\PGLn$-equivariant dominant rational
map such that $h^* = \alpha$.
Then $(f \circ h \circ g\inv)^* = (g\inv)^* \circ \alpha \circ f^* = \beta$.
By uniqueness of $\beta_*$, we have
$f \circ h \circ g\inv = \beta_*$, i.e.,
$h = f\inv \circ \beta_* \circ g = \alpha_*$.
This completes the proof
of Lemma~\ref{lem.prop.functor} and thus of
Proposition~\ref{prop.functor}.
\end{proof}

\begin{cor} \label{cor.rat2}
  Let $X$, $Y$ and $Z$ be generically free irreducible $\PGLn$-varieties.
  \begin{itemize}
  \item[\rm(a)] If $f \colon X \dasharrow Y$, $g \colon Y \dasharrow Z$
  are $\PGLn$-equivariant dominant rational maps then
  $(g \circ f)^* = f^* \circ g^*$.
  \item[\rm(b)]If $\alpha \colon k_n(Y) \hookrightarrow k_n(X)$ and
  $\beta \colon k_n(Z) \hookrightarrow k_n(Y)$ are homomorphisms
  {\upshape(}i.e., embeddings{\upshape)} of central simple algebras
  then $(\alpha \circ \beta)_* = \beta_* \circ \alpha_*$.
  \item[\rm(c)]$X$ and $Y$ are birationally isomorphic as
  $\PGLn$-varieties if and only if $k_n(X)$
  and $k_n(Y)$ are isomorphic as $k$-algebras.
  \end{itemize}
\end{cor}

\begin{proof} (a) is immediate from Definition~\ref{def.gener}.

(b) Let $f = (\alpha \circ \beta)_*$ and
$g  = \beta_* \circ \alpha_*$. Part (a) tells us that
$f^* = g^*$. The uniqueness assertion of
Proposition~\ref{prop.functor} now implies $f = g$.

(c) follows from (a) and (b) and the identities
$(\id_X)^* = \id_{k_n(X)}$, and $(\id_{k_n(X)})_* = \id_X$.
\end{proof}

\begin{proof}[Proof of Theorem~\ref{thm2}]
  We use the same argument as in the proof of Theorem~\ref{thm1}.
  The contravariant functor~$\mathcal F$ is well defined
  by Corollary~\ref{cor.rat2}. Since
  $(f^*)_*=f$ and $(\alpha_*)^*=\alpha$ (see
  Proposition~\ref{prop.functor}), $\mathcal F$ is full and faithful.
  By Theorem~\ref{thm.W(A)}, every object in $\CS_n$ is isomorphic
  to the image of an object in $\Bir_n$.  The desired conclusion
  now follows from~\cite[Theorem~7.6]{blyth}.
\end{proof}

\section{Brauer-Severi Varieties}
\label{sect.bs}

Let $K/k$ be a finitely generated field extension.
Recall that the following sets are in a natural
(i.e., functorial in $K$) bijective correspondence:

\begin{itemize}
\item[(1)] the Galois cohomology set $H^1(K, \PGLn)$,
\item[(2)] central simple algebras $A$ of degree $n$ with center $K$,
\item[(3)] Brauer-Severi varieties over $K$ of dimension~$n-1$,
\item[(4)] $\PGLn$-torsors over $\Spec(K)$,
\item[(5)] pairs $(X, \phi)$, where is $X$ is an irreducible
  generically free $\PGLn$-variety and $\phi_X \colon k(X)^{\PGLn}
  \stackrel{\simeq}{\lra} K$ is an isomorphism of fields (over $k$).
  Two such pairs $(X, \phi)$ and $(Y, \psi)$ are equivalent, if there
  is a $\PGLn$-equivariant birational isomorphism $f \colon Y
  \dasharrow X$ which is compatible with $\phi$ and $\psi$, i.e.,
  there is a commutative diagram
  \[\xymatrix{k(X)^{\PGLn} \ar@{->}[dr]^{\phi}_{\simeq} 
    \ar@{->}[rr]^{f^*}  & &
    k(Y)^{\PGLn} \ar@{->}[dl]^{\simeq}_{\psi} \cr
  & K   & }\]
\end{itemize}
Bijective correspondences between (1), (2), (3) and (4) follow
from the theory of descent; see~\cite[Sections I.5 and III.1]{serre-gc},
\cite[Chapter X]{serre-lf}, \cite[(1.4)]{artin3}
or~\cite[Sections 28, 29]{boi}. For a bijective correspondence
between (1) and (5), see \cite[(1.3)]{popov}.
                             
For notational simplicity we will talk of generically free
$\PGLn$-varieties $X$ instead of pairs $(X, \phi)$ in (5),
and we will write $k(X)^{\PGLn} = K$ instead 
of $k(X)^{\PGLn} \stackrel{\phi}{\simeq} K$, 
keeping $\phi$ in the background.

Suppose we are given a generically free $\PGLn$-variety $X$ (as in (5)).
Then this variety defines a class $\alpha \in H^1(K, \PGLn)$ (as in (1))
and using this class we can recover all the other associated 
objects (2) - (4). In the previous section we saw that
the central simple algebra $A$ can be constructed directly 
from $X$, as $\RMaps_{\PGLn}(X, \Mn)$. The goal of this section
is to describe a way to pass directly from (5) to (3), without 
going through (1); this is done in Proposition~\ref{prop.bs} below.

In order to state Proposition~\ref{prop.bs}, we introduce some 
notation. We will write points of the projective space
$\bbP^{n-1} = \bbP_k^{n-1}$ as rows
$a = (a_1: \ldots : a_n)$. The group $\PGLn$ acts on $\bbP^{n-1}$
by multiplication on the right:
\[ g \colon (a_1: \ldots: a_n) \mapsto (a_1: \ldots : a_n) g^{-1} \, . \]
Choose (and fix) $a = (a_1: \ldots : a_n) \in \bbP^{n-1}$
and define the maximal parabolic subgroup $H$ of $\PGLn$ by
\begin{equation} \label{e.H}
H = \{ h \in \PGLn \, | \,  a h^{-1} = a  \} \, . 
\end{equation}
If $a = (1:0: \dots : 0)$ then $H \subset \PGLn$ consists of
of $n \times n$-matrices of the form 
\[ \begin{pmatrix} * & 0 & \dots & 0 \\
* & * & \dots & * \\
\vdots & \vdots & & \vdots \\
* & * & \dots & * \end{pmatrix} \]

\begin{prop} \label{prop.bs} Let $X$ be an irreducible generically free
$\PGLn$-variety, $A = k_n(X)$ and $\sigma \colon X/H \dasharrow X/\PGLn$.
Then the Brauer-Severi variety $\BS(A)$ is the preimage of the generic point 
$\eta$ of $X/\PGLn$ under $\sigma$. 
\end{prop}

Before we proceed with the proof, three remarks are in order. First
of all, by $X/H$ we mean the rational quotient variety for the 
$H$-action on $X$. Recall that $X/H$ is defined (up to birational 
isomorphism) by $k(X/H) = k(X)^H$, and the dominant rational 
map $\sigma \colon X/H \dasharrow X/\PGLn$ by the 
inclusion of fields $k(X)^{\PGLn} \hookrightarrow k(X)^H$.
Secondly, recall that $k(X/\PGLn) = k(X)^{\PGLn} = K$, so that $\eta
\simeq \Spec(K)$, and $\sigma^{-1}(\eta)$ is, indeed, a $K$-variety. 
Thirdly, while the construction of $\BS(A)$ in Proposition~\ref{prop.bs} 
does not use the Galois cohomology set $H^1(K, \PGLn)$, our proof below 
does.  In fact, our argument is based on showing that 
$\sigma^{-1}(\eta)$ and $\BS(A)$ are Brauer-Severi varieties
defined by the same class in $H^1(K, \PGLn)$. 

\begin{proof} 
Let $X_0/k$ be an algebraic variety with function field $k(X_0) = K$,
i.e., a particular model for the rational quotient variety $X/\PGLn$.
The inclusion $K \stackrel{\phi}{\simeq} k(X)^{\PGLn} \hookrightarrow X$
induces the rational quotient map $\pi \colon X \dasharrow X_0$. 
After replacing $X_0$ by a Zariski dense open subset, we may 
assume $\pi$ is regular; after passing to another (smaller) dense 
open subset, we may assume $\pi \colon X \lra X_0$ is,
in fact, a torsor; cf., e.g., \cite[Section 2.5]{pv}.           
                             
We now trivialize this torsor over some etale cover $U_i \to X_0$.
Then for each $i, j$ the transition map 
$f_{ij} \colon \PGLn \times U_{ij} \lra \PGLn \times U_{ij}$
is an automorphism of the trivial
$\PGLn$-torsor $\PGLn \times U_{ij}$ on $U_{ij}$. 
It is easy to see that $f_{ij}$ is given by the formula
\begin{equation} \label{e.h} 
f_{ij} (u, g) = (u, g \cdot c_{ij}(u)) \, ,
\end{equation}               
for some morphism $c_{ij} \colon U_{ij} \lra \PGLn$.
The morphisms $c_{ij}$ satisfy a cocycle condition
(for Cech cohomology) which expresses the fact that the
transition maps $f_{ij}$ are compatible on triple
``overlaps" $U_{hij}$. The cocycle
$c = (c_{ij})$ gives rise to a cohomology
class $\overline{c} \in H^1(X_0, \PGLn)$, which maps to $\alpha$
under the natural restriction morphism $H^1(X_0, \PGLn) \lra H^1(K, \PGLn)$
from $X_0$ to its generic point; 
cf.\ \cite[Section 8]{cgr}.  (Recall that by our construction 
the function field of $X_0$ is identified with $K$.)
                             
Now define the quotient $Z$ of $X$ by the maximal
parabolic subgroup $H \subset \PGLn$ as follows.
Over each $U_i$ set $Z_i = H \backslash \PGLn \times U_i$. By descent
we can ``glue" the projection morphisms $Z_i \lra U_i$ into
a morphism $Z \lra X_0$ by the transition maps
\[ \overline{f_{ij}} (u, g) = (u, \overline{g} \cdot c_{ij}(u)) \, . \]
Moreover, since over each $U_i$ the map $\pi \colon X \lra X_0$ factors as
\[ \pi_i \colon \PGLn \times U_i \stackrel{p_i}{\lra}
H \backslash \PGLn \times U_i \stackrel{q_i}{\lra} U_i \, , \]
the projection maps $p_i$ and $q_i$ also glue together,
yielding                     
\[ \pi \colon X \stackrel{p}{\lra} Z \stackrel{q}{\lra} X_0 \, . \]

By our construction the fibers of $p$ are exactly
the $H$-orbits in $X$; hence, $k(Z) = k(X)^H$, cf., e.g.,~\cite[2.1]{pv}. 
In other words, $p$ is a rational quotient map for 
the $H$-action on $X$ and we can identify $Z$ with the rational 
quotient variety $X/H$ (up to birational equivalence).
Under this identification $q$ becomes $\sigma$.
                             
Now recall that by the definition of $H$, the homogeneous 
space $H \backslash \PGLn$ is naturally
isomorphic with $\bbP^{n-1}$ via $\overline{g} \mapsto g \cdot a$.
Since over each $U_i$ the map $q \colon Z \lra X_0$ looks like
the projection $H \backslash \PGLn \times U_i \lra U_i$,
$Z$ is, by definition, a Brauer-Severi variety over $X_0$.
Moreover, $\pi \colon X \lra X_0$ (viewed as a torsor over $X_0$) and 
$q \colon Z \lra X_0$ (viewed as a Brauer-Severi variety over $X_0$)
are constructed by using the same cocycle $(c_{ij})$
and hence, the same cohomology class $\overline{c} \in H^1(X_0, \PGLn)$.
Restricting to the generic point of $X_0$, we see that the cohomology
class of $Z$ as a Brauer-Severi variety over $K = k(X_0)$ is the image of
$\overline{c}$ under the restriction map $H^1(X_0, \PGLn) \lra H^1(K, \PGLn)$, 
i.e., the class $\alpha \in H^1(K, \PGLn)$ we started out with. 
\end{proof}                   
         
\begin{remark} Note that the choice of the maximal parabolic
subgroup $H \subset \PGLn$ is important here. If we repeat
the same construction with $H$ replaced by $H^{\text{transpose}}$
we will obtain the Brauer-Severi variety
of the opposite algebra $A^{\text{op}}$.
\end{remark}                 

The following corollary of Proposition~\ref{prop.bs} shows that $X$,
viewed as an abstract variety (i.e., without the $\PGLn$-action), is
closely related to $\BS(A)$.

\begin{cor} \label{cor.bs} Let $X$ be a generically free $\PGLn$-variety,
$A = k_n(X)$ be the associated central simple algebra of degree $n$,
and $K = k(X)^{\PGLn}$ be the center of $A$. Then 
\[ k(X) \simeq K(\BS(A))(t_1, \dots, t_{n^2-n}) \simeq K(\BS(\Mn(A))) \, . \] 
\end{cor}

Here $k(X)$ denotes the function field of $X$ (as a
variety over $k$), and $K(\BS(A))$ and $K(\BS(\Mn(A))$ 
denote, respectively, the function fields of the 
Brauer-Severi varieties of $A$ and $\Mn(A)$ (both are
defined over $K$). The letters $t_1, \dots, t_{n^2-1}$ denote $n^2 - n$ 
independent commuting variables, and the isomorphisms $\simeq$ are field
isomorphisms over $k$ (they ignore the $\PGLn$-action on $X$).

\begin{proof} The second isomorphism is due to 
Roquette~\cite[Theorem 4, p. 413]{roquette}. 
To show that
$k(X) \simeq K(\BS(A))(t_1, \dots, t_{n^2-n})$,
note that by Proposition~\ref{prop.bs}, $K(\BS(A)) = K(X/H) = k(X)^H$, 
where $H$ is the parabolic subgroup of $\PGLn$
defined in~\eqref{e.H}. Since $\dim(H) = n^2 - n$, it
remains to show that the field extension $k(X)/k(X)^H$ 
is rational. 

Now recall that $H$ is a special group (cf.~\cite[Section 2.6]{pv}); 
indeed, the Levi subgroup of $H$ is isomorphic to
$\GL_{n-1}$. Consequently, $X$ is birationally isomorphic 
to $(X/H) \times H$ (over $k$). Since $k$ is assumed to be 
algebraically closed and of characteristic zero, every 
algebraic group over $k$ is rational. In particular, $H$ 
is birationally isomorphic to
$\bbA^{n^2-n}$ and thus $X$ is birationally isomorphic to 
$(X/H) \times \bbA^{n^2-n}$.
In other words,
\[ k(X) \simeq K(X/H)(t_1, \dots, t_{n^2-n}) \simeq 
K(\BS(A))(t_1, \dots, t_{n^2-n}) \, , \]
as claimed.
\end{proof}


\begin{thebibliography}{KMRT}

\bibitem{amitsur}S. A. Amitsur, \emph{A generalization of
    Hilbert's Nullstellensatz}, Proc. Amer. Math. Soc. {\bf 8}
  (1957), 649--656.
\myendbibitem

\bibitem{amitsur:procesi}S. A. Amitsur and C. Procesi,
  \emph{Jacobson-rings and Hilbert algebras with polynomial
  identities},
Ann. Mat. Pura Appl. (4) {\bf 71} (1966), 61--72.
\myendbibitem

\bibitem{amitsur:small}S. A. Amitsur and L. W. Small,
  \emph{Prime ideals in PI-rings},
J. Algebra {\bf 62} (1980), 358--383.
\myendbibitem

\bibitem{artin1}M. Artin, \emph{On Azumaya algebras
and finite dimensional representations of rings},
J. Algebra {\bf 11} (1969), 532--563.
\myendbibitem

\bibitem{artin2} M. Artin,  {\em Specialization of
representations of rings}, Proceedings of the International
Symposium on Algebraic Geometry (Kyoto Univ., Kyoto, 1977),
Kinokuniya Book Store, Tokyo (1978), 237-247.
\myendbibitem
                             
\bibitem{artin3} M.  Artin, \emph{Brauer-Severi varieties},
Brauer groups in ring theory and algebraic geometry (Wilrijk, 1981),
pp. 194--210, Lecture Notes in Math., 917,
Springer, Berlin-New York, 1982.
\myendbibitem

\bibitem{as1} M. Artin, W. Schelter,
{\em A version of Zariski's main theorem for polynomial identity rings}.
Amer. J. Math. 101 (1979), no. 2, 301--330.
\myendbibitem

\bibitem{blyth}T. S. Blyth, {\em Categories}, Longman, London and
  New York, 1986.
\myendbibitem

\bibitem{cgr} V. Chernousov, Ph. Gille, Z. Reichstein,
Resolving G-torsors by abelian base extensions, ArXIV math.AG/0404392,
J. Algebra, to appear. 
\myendbibitem

\bibitem{demeyer}F. DeMeyer, E. Ingraham, \emph{Separable Algebras
    over Commutative Rings}, Lecture Notes in Mathematics {\bf 181},
    Springer-Verlag, 1971.
\myendbibitem

\bibitem{formanek:cbms}E. Formanek, \emph{The polynomial identities
    and invariants of $n\times n$ matrices}, CBMS Regional Conference
    Series in Mathematics {\bf 78}, 1991.
\myendbibitem

\bibitem{friedland}S. Friedland, 
\emph{Simultaneous similarity of matrices},
Adv. in Math. {\bf 50}  (1983), no. 3, 189--265.
\myendbibitem

\bibitem{hartshorne}R. Hartshorne, \emph{Algebraic Geometry},
Springer-Verlag, 1977.
\myendbibitem

\bibitem{boi} M.-A.
 Knus, A. Merkurjev, M. Rost, J.-P. Tignol,
\emph{The book of involutions}, Colloquium Publications, {\bf 44},
American Mathematical Society, Providence, RI, 1998.
\myendbibitem

\bibitem{lebruyn}L. Le Bruyn, {\em Local structure 
of Schelter-Procesi smooth orders},
Trans. Amer. Math. Soc. 352 (2000), no. 10, 4815--4841.
\myendbibitem

\bibitem{mf}D. Mumford, J. Fogarty, {\em Geometric Invariant Theory},
second edition, Springer-Verlag, 1982.
\myendbibitem

\bibitem{popov}
V. L. Popov, \emph{Sections in invariant theory}, The
Sophus Lie Memorial Conference (Oslo, 1992), Scand. Univ. Press,
Oslo, 1994, pp.~315--361.
\myendbibitem

\bibitem{pv}V. L. Popov and E. B. Vinberg, {\it Invariant Theory},
Algebraic Geometry IV, Encyclopedia of Mathematical Sciences {\bf
55}, Springer, 1994, 123--284.
\myendbibitem

\bibitem{procesi2}C. Procesi, {\em Rings with polynomial
identities}, Marcel Dekker Inc., 1973.
\myendbibitem

\bibitem{procesi3}C. Procesi, {\em Finite dimensional 
representations of algebras}, Israel J. Math. 19 (1974), 169--182.
\myendbibitem

\bibitem{procesi1}C. Procesi, \emph{The
invariant theory of $n\times n$ matrices}, Advances in Math. {\bf 19}
(1976), no. 3, 306--381.
\myendbibitem

\bibitem{procesi4}C. Procesi, {\em A formal inverse 
to the Cayley-Hamilton theorem}, J. Algebra 107 (1987), no. 1, 63--74.
\myendbibitem

\bibitem{razmyslov} Yu. P.  Razmyslov, 
The Jacobson Radical in PI-algebras. (Russian)
Algebra i Logika {\bf 13} (1974), 337--360, 365.
English translation: Algebra and Logic {\bf 13} (1974), 192--204 (1975).
\myendbibitem

\bibitem{reichstein4} Z.~Reichstein, \emph{On the notion of essential
    dimension for algebraic groups}, Transform. Groups {\bf 5}
  (2000), no.~3, 265--304.
\myendbibitem

\bibitem{roquette} P. Roquette, {\em On the Galois cohomology 
of the projective linear group and its applications to the construction 
of generic splitting fields of algebras}
Math. Ann. {\bf 150} (1963), 411--439.
\myendbibitem

\bibitem{rowen:PI}L.\ H.\ Rowen, {\it Polynomial Identities in Ring
    Theory}, Academic Press, 1980.
\myendbibitem

\bibitem{rowen:ringII}L.\ H.\ Rowen, {\it Ring Theory}, Volume II,
  Academic Press, 1988.
\myendbibitem

\bibitem{schelter:78}W. Schelter, \emph{Non-Commtutative Affine
    P. I. Rings are Catenary}, J. Algebra {\bf51} (1978), 12--18.
\myendbibitem

\bibitem{serre-gc}J - P.~Serre, {\it Galois Cohomology},
Springer, 1997.
\myendbibitem

\bibitem{serre-lf}J.-P. Serre, {\it Local Fields}, 
Springer -- Verlag, 1979.
\myendbibitem

\end{thebibliography}
\end{document}